\documentclass[]{interact}

\usepackage[caption=false]{subfig}

\usepackage{amsmath}
\usepackage{amsfonts}
\usepackage{bbm}
\usepackage[american]{babel}
\usepackage[utf8]{inputenc}
\usepackage[T1]{fontenc}
\usepackage[linesnumbered,lined,boxed,commentsnumbered]{algorithm2e}

\usepackage{eurosym}

\usepackage{csquotes}
\usepackage[style=apa, backend=biber]{biblatex}
\usepackage[colorlinks,
            citecolor=blue,
            urlcolor=blue,
            linkcolor=blue]{hyperref} 

\usepackage{tikz}
\usetikzlibrary{positioning}
\usetikzlibrary{calc}
\usetikzlibrary{arrows}
\usetikzlibrary{decorations.pathmorphing,decorations.markings}
\usetikzlibrary{shapes}
\usetikzlibrary{patterns}

\usepackage{pgf-pie}  

\theoremstyle{plain}
\newtheorem{theorem}{Theorem}[section]

\newtheorem{definition}[theorem]{Definition}

\newcommand{\st}{\text{s.t.}}

\newcommand{\set}[1]{\{#1\}}

\DeclareMathOperator*{\argmin}{arg\,min}

\begin{document}

\title{A novel strong duality-based reformulation for trilevel infrastructure models in energy systems development}
\author{
\name{Olli Herrala\textsuperscript{a}, Steven A. Gabriel\textsuperscript{a,b,c}, Fabricio Oliveira\textsuperscript{a,*}, and Tommi Ekholm\textsuperscript{a,d}}
\affil{\textsuperscript{a}Aalto University, Department of Mathematics and Systems Analysis, Espoo, Finland; \textsuperscript{b}University of Maryland, Department of Mechanical Engineering/Applied Mathematics, Statistics \& Scientific Computing Program, College Park, Maryland USA; \textsuperscript{c}Norwegian University of Science and Technology, Department of Industrial Economics and Technology Management, Trondheim, Norway; \textsuperscript{d}Finnish Meteorological Institute, Helsinki, Finland; \textsuperscript{*}Corresponding author: fabricio.oliveira@aalto.fi
}}

\maketitle

\begin{abstract}
We explore the class of trilevel equilibrium problems with a focus on energy-environmental applications and present a novel single-level reformulation for such problems, based on strong duality. To the best of our knowledge, only one alternative single-level reformulation for trilevel problems exists. This reformulation uses a representation of the bottom-level solution set, whereas we propose a reformulation based on strong duality. Our novel reformulation is compared to this existing formulation, discussing both model sizes and computational performance.  In particular, we apply this trilevel framework to a power market model, exploring the possibilities of an international policymaker in reducing emissions of the system. Using the proposed methods, we are able to obtain globally optimal solutions for a five-node case study representing the Nordic countries and assess the impact of a carbon tax on the electricity production portfolio.
\end{abstract}


\section{Introduction}
\label{sec:introduction}

Hierarchical optimisation models with three levels of decision-makers arise in contexts such as traffic equilibrium \parencite{gu2019tri, Gabriel_et_al:2021} and electricity market modelling \parencite{jin2013tri, huppmann2015national}. The hierarchical structure can be, e.g., such that the bottom-level players use a network operated by a middle-level player and regulated by a top-level player. For both electricity and traffic networks, similar models without the top-level regulators have been explored using \emph{bilevel optimisation}, see \textcite{sinha2017review} for a review. 

Albeit challenging from both methodological and computational standpoints, including a top-level regulator as the third level, as opposed to considering only bilevel models, can provide important policy insights. In the particular case of energy systems, these models can yield more realistic solutions in which more stakeholders are assumed to act in coordination considering their own objectives. Obtaining equilibrium solutions for these models can thereby provide policy insights on pathways towards decarbonisation goals. \textcite{Gabriel_et_al:2021} present a single-level reformulation for bilevel problems with complementarity-constrained bottom levels and discuss the possibility of using the model in a trilevel power market setting. However, their paper includes no computational experiments demonstrating the practical usability of the proposed methods. Our aim is to explore the computational efficiency of the method using an illustrative power system setting representing the market structure in the Nordic countries.

The contribution of this paper is twofold. First, in Section \ref{sec:background}, we present background on bi- and tri-level optimisation, ending with our novel approach for solving trilevel equilibrium problems based on strong duality in Section \ref{sec:mpc-pdc}. Compared to the formulation presented in \textcite{Gabriel_et_al:2021} and summarised in Section \ref{sec:posit-semid-m-x-dep}, our proposed formulation results in fewer constraints, which is likely to result in increased computational efficiency. 
Second, we illustrate the methodological contributions using a stylised trilevel power market model described in Section \ref{sec:applications}. The motivation for our application stems from the recent discussion about optimal carbon taxation and its impact on electricity production \parencite[e.g.,][]{hajek2019analysis}.  The computational performance of the model is explored in Section \ref{sec:formulations-comp}, and finally, in Section \ref{sec:casestudy}, we apply the trilevel equilibrium modelling framework to a power market case study based on \textcite{belyak2023optimal}. These contributions are significant for the novel class of trilevel optimisation problems, and equilibrium modelling area in general. Finally, Section \ref{sec:conclusions} concludes the paper and discusses future research directions.

\section{Background}
\label{sec:background}

\subsection{Earlier research}
\label{sec:litrev}

Bilevel optimisation considers problems with a hierarchical structure consisting of an upper-level player and one or more lower-level players \parencite{bard1983algorithm}. In power sector models \parencite{gabriel2012complementarity}, the upper-level player is often a transmission system operator and the bottom level consists of electricity producers in a Cournot oligopoly. The general structure of a bilevel problem with linear upper- and lower-level problems is presented in \eqref{eq:bilevel-upper} and \eqref{eq:bilevel-lower}. The upper-level problem $P_u$ is
\begin{subequations}
  \label{eq:bilevel-upper}
  \begin{align}
    (P_u): \ min_{x,y \ge 0} \quad & c_1^\top x + d_1^\top y \label{eq:bilevel-upper.objective} \\
    \st \quad & A_1 x + B_1 y \geq a_1, \label{eq:bilevel-upper.constraints} \\
    & y \text{ solves } P_l(x) \label{eq:bilevel-upper.lowerlevel},
  \end{align}
\end{subequations}
where $P_l(x)$ denotes the lower-level problem. Here, $c_1, x \in R^{n_x}$, $d_1, y \in R^{n_y}$, $A_1 \in R^{m_u \times n_x}$, $B_1 \in R^{m_u \times n_y}$, and $a_1 \in R^{m_u}$. The overall idea of this formulation is that the upper-level player's decision variable $x$ affects the lower-level players' optimal decisions $y$, which are reflected back to the upper level in the constraint \eqref{eq:bilevel-upper.lowerlevel}. The linear lower-level problem is formulated as
\begin{subequations}
  \label{eq:bilevel-lower}
  \begin{align}
    (P_l(x)): \ min_{y \ge 0} \quad & d_2^\top y \label{eq:bilevel-lower.obj} \\
    \st \quad & A_2 x + B_2 y \geq a_2. \label{eq:bilevel-lower.cons}
  \end{align}
\end{subequations}
In general, both problems can also include equality constraints, but they have been omitted here for brevity, without loss of generality. 

Solution methods for bilevel problems are based on the idea of replacing the upper-level constraint \eqref{eq:bilevel-upper.lowerlevel} with the optimality conditions of the lower level problem \eqref{eq:bilevel-lower}. The two main alternatives are the Karush-Kuhn-Tucker (KKT) optimality conditions \parencite{karush1939minima, kuhn1951nonlinear}, leading to a mathematical program with equilibrium constraints (MPEC); and mathematical programming with primal and dual constraints (MPPDC) \parencite{ruiz2011equilibria, baringo2012transmission}. Additionally, approaches based on optimal value functions \parencite{ye1995optimality} can be used. Bilevel optimisation models can be used in contexts such as Stackelberg games \parencite{bard1991some}, Cournot competition \parencite{gabriel2012complementarity} and robust optimisation \parencite{leyffer2020survey}. For a recent survey on applications and algorithms for bilevel optimisation, we refer to \textcite{kleinert2021survey}.

A major challenge in bilevel optimisation is that, even for linear bilevel problems, single-level reformulations using complementarity constraints lead to nonlinear and nonconvex problems. This significantly increases the computational complexity of solving such problems and requires specialised approaches such as the simplex method-inspired projected gradient method by \textcite{still2002linear} or the (spatial) branch-and-bound methods discussed by \textcite{bard1990branch} and implemented in the Gurobi solver \parencite{gurobi}. Alternatively, one can use heuristics such as genetic programming \parencite{kieffer2019tackling} and particle swarm optimisation \parencite{gao2021improving}. For a review on heuristic solution methods for bilevel programming, we refer the reader to \textcite{camacho2023metaheuristics}.




\subsection{Trilevel equilibrium models}
\label{sec:problem-statement}

Consider a problem with a trilevel structure, in which players interact with each other at all three levels: top, middle and bottom. In this structure, the top-level problem $P_1$ is assumed to be a linear optimisation problem with the middle-level problem $P_2(x)$ represented by the constraint \eqref{eq:toplevel-prob.middlelevel}.
\begin{subequations}
  \label{eq:toplevel-prob}
  \begin{align}
    (P_1): \ min_{x,y,z} \quad & c_1^\top x + d_1^\top y + e_1^\top z \\
    \st \quad & A_1 x + B_1 y + C_1 z \geq a_1 \label{eq:toplevel-prob.constraints}\\
    & y,z \text{ solve } P_2(x), \label{eq:toplevel-prob.middlelevel}
  \end{align}
\end{subequations}
where $P_2(x)$ denotes the middle-level problem
\begin{subequations}
  \label{eq:middlelevel-prob}
  \begin{align}
    (P_2(x)): \ min_{y,z} \quad & d_2^\top y + e_2^\top z \label{eq:middlelevel-prob.obj} \\
    \st \quad & A_2 x + B_2 y + C_2 z \geq a_2 \quad (\gamma) \label{eq:middlelevel-prob.cons} \\
    & y \ge 0 \label{eq:middlelevel-prob.nonneg} \\
    & z \text{ solves } P_3(x,y), \label{eq:middlelevel-prob.bottomlevel}
  \end{align}
\end{subequations}
where $\gamma$ is the vector of dual variables associated with constraint \eqref{eq:middlelevel-prob.cons} and $z$ is the vector of bottom-level variables. In trilevel settings, the ``lower-level'' problem $P_2(x)$ is itself a bilevel problem. This is challenging, because bilevel optimisation problems are generally nonconvex and directly obtaining their optimality conditions is thus difficult. The middle-level problem $P_2(x)$ constraints contain a bottom-level problem $P_3(x,y)$ that is parameterised by the upper-level variables~$x$ and middle-level variables $y$. In particular, \textcite{Gabriel_et_al:2021} discuss the case of the bottom-level problem being a linear complementarity problem (LCP) 
\begin{equation} 
P_3(x,y): \quad 0 \leq \tilde{z} \perp q + N_x x + N_y y + M\tilde{z} \geq 0 \label{eq:lower-level-lcp}
\end{equation}
parameterised via the vector terms $N_x x$ and $N_y y$. It should be noted that \eqref{eq:lower-level-lcp} can be viewed as the KKT conditions of convex quadratic problems and in Section \ref{sec:mpc-pdc}, we discuss such problems in more detail. Hereinafter, we use the standard $\perp$-notation
\begin{equation*}
  0 \leq a \perp b \geq 0
  \iff
  a, b \geq 0, \, a^\top b = 0
\end{equation*}
for complementarity constraints with vectors $a$ and $b$.



Trilevel problems have been researched by, e.g., \textcite{sauma2007economic}, who do not solve the models directly and instead iteratively solve the middle- and bottom-level problem for different values of the top-level decision variables; and \textcite{dvorkin2017co}, who employ a column-and-constraint generation algorithm. In contrast, the goal of this paper is to explore novel single-level reformulations for trilevel problems and solve them using an off-the-shelf solver, thus lowering the barrier-to-entry for using these models. However, this comes at the expense of imposing constraints on the structure of the models that are tractable in this manner. In \eqref{eq:middlelevel-prob}, a general form of the problem is used, and the bottom-level problem in \eqref{eq:middlelevel-prob.bottomlevel} is assumed to be parameterised by both $x$ and $y$. For the sake of clarity, we define two classes of trilevel problems with different degrees of computational challenges.

\begin{definition}
  \label{def:structure}
    If $P_3$ is parameterised by both $x$ and $y$, we say that the problem has a \emph{strong} trilevel structure. In contrast, if $P_3$ is parameterised only by the top-level variables $x$, i.e., it is not directly dependent on $y$, we say that the problem has a \emph{weak} trilevel structure.
\end{definition}

\textcite{Gabriel_et_al:2021} show that in order to use their single-level reformulation, the problem must have a weak trilevel structure, allowing such problems to be solved rather effectively by borrowing from the results in \textcite[Theorem~3.1.6]{Cottle_et_al:2009} as long as the matrix $M$ in the lower-level problem \eqref{eq:lower-level-lcp} is positive semi-definite (PSD). We also show that our novel reformulation in Section \ref{sec:mpc-pdc} retains this structural limitation, and that the energy-environmental planning problem considered in this paper has this structure. The aim of this paper is to develop an alternative reformulation improving computational tractability and efficiency compared to the reformulation in \textcite{Gabriel_et_al:2021}, and the discussion on ways for lifting this restriction on problem structure is outside the scope of this paper. 

While the lack of direct influence for the middle-level player is a limitation, there are still structures that necessitate the use of a trilevel framework. 
As an example of a setting where a trilevel approach is required, we use the power market example in Section \ref{sec:applications}, where the bottom level consists of electricity generators, and on the middle level we have a profit-maximising system operator who has to satisfy a minimum renewable share in electricity production. If the bottom-level LCP matrix $M$ in \eqref{eq:lower-level-lcp} is PSD, the bottom-level problem can have multiple optima. This could result in, e.g., a situation where it makes no difference for a generator to produce electricity using coal in one node or wind power in another. 

Using the optimistic bilevel assumption \parencite{dempe2020bilevel}, while the system operator cannot directly influence the generators, they can choose a bottom-level optimum that maximises their profit while satisfying the minimum renewable share constraint. In turn, maximising the middle-level player's profit could, in some settings, result in worse objective values for the top-level player. These interactions could not be represented in a setting where the middle-level player is insensitive to the bottom-level player's decisions, as the middle-level player must consider the bottom-level optimality conditions to be able to choose between bottom-level optimal solutions.

\subsection{Bottom-level LCP with a positive semi-definite coefficient matrix}
\label{sec:posit-semid-m-x-dep}

For completeness, we summarise the solution approach introduced in \textcite{Gabriel_et_al:2021}. Let us assume that the matrix $M$ in \eqref{eq:lower-level-lcp} is positive semi-definite and that we have a solution $\bar{z}$ of \eqref{eq:lower-level-lcp}. Furthermore, we assume the problem to have a weak trilevel structure and thus $N_y=0$, i.e., the middle-level decisions $y$ do not influence the bottom-level problem (c.f. Definition \ref{def:structure}). 
\textcite{Gabriel_et_al:2021} show that for a positive semi-definite $M$ and a weak trilevel structure, a solution to the trilevel problem consisting of \eqref{eq:toplevel-prob}-\eqref{eq:lower-level-lcp} can be obtained by solving the equivalent single-level reformulation
\begin{subequations}
\label{eq:large-single-level-lpcc-in-x-and-y}
    \begin{align}
    \min_{x, y, \tilde{z}, \bar{z}, \beta, \gamma, \delta, \zeta, \eta} \quad & c_1^\top x + d_1^\top y + e_1^\top \tilde{z} \label{subeq:large-single-level-lpcc-in-x-and-y.a} \\
    \st \quad & A_1 x + B_1 y + C_1 \tilde{z} \geq a_1, \label{subeq:large-single-level-lpcc-in-x-and-y.b} \\
    & 0 \leq y \perp d_2 - B_2^\top \gamma \geq 0, \label{subeq:large-single-level-lpcc-in-x-and-y.g1} \\
    & 0 \leq \tilde{z} \perp e_2 - C_2^\top \gamma - M^\top \delta - \zeta (q + N_x x) - (M + M^\top)^\top \eta \geq 0, \label{subeq:large-single-level-lpcc-in-x-and-y.g2} \\
    & 0 \leq \bar{z} \perp q + N_x x + (M+M^\top)\bar{z} - M^\top \beta \geq 0, \label{subeq:large-single-level-lpcc-in-x-and-y.c} \\
    & 0 \leq \beta \perp q + N_x x + M \bar{z} \geq 0, \label{subeq:large-single-level-lpcc-in-x-and-y.d} \\
    & 0 \leq \delta \perp (q + N_x x) + M\tilde{z} \geq 0, \label{subeq:large-single-level-lpcc-in-x-and-y.h} \\
    & (q + N_x x)^\top(\tilde{z} - \bar{z}) = 0, \ (M + M^\top)(\tilde{z} - \bar{z})  = 0, \label{subeq:large-single-level-lpcc-in-x-and-y.i} \\
    & 0 \leq \gamma \perp A_2 x + B_2 y + C_2 \tilde{z} - a_2 \geq 0, \label{subeq:large-single-level-lpcc-in-x-and-y.j}
    \end{align}
\end{subequations}
where $\bar{z}$ is a solution to the bottom-level problem \eqref{eq:lower-level-lcp} and $(\bar{z}^*)^\top (q + N_x x^* + M\bar{z}^*) = 0$ has to thus hold at an optimal solution $x^*, \bar{z}^*$ to \eqref{eq:large-single-level-lpcc-in-x-and-y}. Appendix \ref{app:posit-semid-m-x-dep} summarizes the reformulation steps taken in \textcite{Gabriel_et_al:2021}, including the constraints corresponding to the dual variables $\beta, \delta, \zeta$ and $\eta$. This formulation assumes nonnegativity for all variables $y$, but we note that this is not a requirement and including free variables in the middle level only requires small changes to the corresponding KKT conditions \eqref{subeq:large-single-level-lpcc-in-x-and-y.g1}.

Finally, we note that the reformulation \eqref{eq:large-single-level-lpcc-in-x-and-y} is not linear due to the nonlinear products $\zeta N_x x$, $x^\top N_x^\top \tilde{z}$ and $x^\top N_x^\top \bar{z}$, resulting in a nonconvex problem. In general, obtaining global optimal solutions to nonconvex problems is enormously challenging, but this particular nonconvexity can be handled by using a solver capable of handling problems with bilinear terms in special ordered sets of type 1 (SOS1) constraints \parencite{beale1970special}. An SOS1 constraint states that out of a set of variables or functions, only one can have a nonzero value. A complementarity constraint $0 \le a \perp b \ge 0$ can thus be reformulated as two nonnegative variables $a$ and $b$ in a SOS1 constraint. This can be achieved using, e.g., the spatial branch-and-bound method in the Gurobi solver \parencite{gurobi}, see also \textcite{siddiqui2013sos1}. 

\subsection{Mathematical programming with complementarity from primal and dual constraints}
\label{sec:mpc-pdc}

We are now ready to discuss our novel single-level reformulation. Let us first consider a setting where the bottom level is a convex quadratic minimization problem. So far, we have discussed a reformulation based on adding the KKT optimality conditions of the bottom-level problem to the middle-level problem. In our trilevel case, the KKT optimality conditions, having complementarity constraints, require a reformulation of the LCP solution set so that we can obtain a single-level equivalent formulation of the trilevel problem. This eventually results in the middle- and bottom-level problems being represented as two optimisation problems, potentially leading to computational challenges with the reformulation \eqref{eq:large-single-level-lpcc-in-x-and-y}. Representing these two nested optimisation problems as a single-level equivalent requires a large number of complementarity constraints \eqref{subeq:large-single-level-lpcc-in-x-and-y.c}-\eqref{subeq:large-single-level-lpcc-in-x-and-y.j}, possibly leading to prohibitive computational requirements. 

To circumvent these challenges, we note that some bilevel optimisation problems can also be reformulated as mathematical programs with primal and dual constraints (MPPDC), using strong duality instead of complementarity. We present a novel strong duality-based reformulation for trilevel problems, in which a linear middle-level problem and convex quadratic bottom-level problems are reformulated into a single quadratically constrained linear problem (QCLP) instead of two optimisation problems (a QP and an LP) as in Appendix \ref{app:posit-semid-m-x-dep} and \textcite{Gabriel_et_al:2021}. The model sizes resulting from using complementarity (Section \ref{sec:posit-semid-m-x-dep}) and strong duality (this section) for the bottom level are compared in Section \ref{sec:comparison}.

Consider a trilevel problem with a set of bottom-level problems $P_{3i}(x)$
\begin{subequations}
\label{eq:mppdc_bottom_primal}
\begin{align}
  (P_{3i}(x)): \ \min_{z_i} \quad & \frac{1}{2} z_i^\top F_i z_i + e_{i3}(x)^\top z_i \\
  \st \quad & C_{i3} z_i \ge a_{i3}(x) \\
  & z_i \geq 0,
\end{align}
\end{subequations}
where $z_i$ is a vector of decision variables and $F_i$ is positive semidefinite (PSD) for all $i \in I$. In our illustrative example described in the next section, the set $I$ represents the electricity producers. Note that we assume a weak trilevel structure, that is, $P_{3i}$ does not depend on $y$. \textcite{dorn1960duality} presents Lagrangian dual formulations for quadratic problems\footnote{For completeness, the steps for obtaining the dual \eqref{eq:mppdc_bottom_dual} from the primal problem \eqref{eq:mppdc_bottom_primal} are presented in Appendix \ref{app:dual-formulation}.}, and using these formulations, the dual of each problem $P_{3i}(x)$ is 
\begin{subequations}
\label{eq:mppdc_bottom_dual}
\begin{align}
  \max_{p_i, z_i} \quad & -\frac{1}{2} z_i^\top F_i z_i + a_{i3}(x)^\top p_i \\
  \st \quad & C_{i3}^\top p_i - F_i z_i \le e_{i3}(x)\\
  & p_i \geq 0.
\end{align}
\end{subequations}
In MPPDC, the complementarity constraints in the KKT optimality conditions are replaced with a \emph{strong duality} constraint. The strong duality theorem \parencite[e.g.,][]{bazaraa2013nonlinear} states that if the problem has no duality gap, that is, some constraint qualification holds for the problem\footnote{For problems with only affine constraints, such as \eqref{eq:mppdc_bottom_primal}, the Abadie constraint qualification is always satisfied \parencite{bazaraa2013nonlinear}.}, the optimal primal and dual objective values are equal. This implies that such problems can be solved to optimality by finding any solution that is both primal and dual feasible with the primal and dual objective values being equal.

Combining formulations \eqref{eq:mppdc_bottom_primal} and \eqref{eq:mppdc_bottom_dual}, we obtain the following primal and dual constraints, combined with a strong-duality constraint:
\begin{subequations}
\label{eq:bottom-level-strong-duality}
\begin{align}
    & C_{i3} z_i \ge a_{i3}(x) \ \forall i \in I \\
    & C_{i3}^\top p_i - F_i z_i \le e_{i3}(x) \ \forall i \in I \\
    & z_i^\top F_i z_i + e_{i3}(x)^\top z_i -  a_{i3}(x)^\top p_i \le 0 \ \forall i \in I \label{eq:strong-duality-separated} \\
    & z_i, p_i \geq 0 \ \forall i \in I.
\end{align}
\end{subequations}
The strong duality constraint \eqref{eq:strong-duality-separated} states that the objective value of each bottom-level primal (minimisation) problem must not be higher than the value of the dual (maximisation) problem. Recall that the weak duality theorem \parencite{bazaraa2013nonlinear} states that the objective value of any solution of a minimisation problem is greater or equal to any objective value of the corresponding dual problem. This result allows us to write the strong duality constraint in an inequality form, following the approach in \textcite{huppmann2015national}, thus avoiding a quadratic equality constraint that would render a nonconvex feasible region. Since the matrices $F_i$ are PSD, constraints \eqref{eq:strong-duality-separated} are convex. Knowing that weak duality guarantees the left-hand side of each constraint \eqref{eq:strong-duality-separated} to be nonnegative also allows us to combine the $|I|$ constraints into one by taking a sum over the left-hand side values, reducing the number of constraints. 

By combining the middle-level problem \eqref{eq:middlelevel-prob.obj}-\eqref{eq:middlelevel-prob.nonneg} with the bottom-level problem reformulation \eqref{eq:bottom-level-strong-duality}, we obtain the resulting bilevel MPPDC formulation of \eqref{eq:middlelevel-prob}:
\begin{subequations}
\label{eq:middlelevel-mppdc}
\begin{align}
    \min_{y, z_i, p_i} \quad & d_2^\top y + \sum_{i \in I} e_{i2}^\top z_i \\
    \st \quad & A_2 x + B_{2} y + \sum_{i \in I} C_{i2} z_i \le a_2 \\
    & C_{i3} z_i \ge a_{i3}(x) \ \forall i \in I \label{eq:middlelevel-mppdc-con2}\\
    & C_{i3}^\top p_i - F_i z_i \le e_{i3}(x) \ \forall i \in I \label{eq:middlelevel-mppdc-con3}\\
    & \sum_{i \in I} (z_i^\top F_i z_i + e_{i3}(x)^\top z_i -  a_{i3}(x)^\top p_i) \le 0 \label{eq:strong-duality-mppdc} \\
    & y \ge 0 \\
    & z_i, p_i \geq 0 \ \forall i \in I.
\end{align}
\end{subequations}
The objective function \eqref{eq:middlelevel-prob.obj} and constraint \eqref{eq:middlelevel-prob.cons} have been modified from \eqref{eq:middlelevel-prob} by adding a sum over the set $I$ to highlight the fact that we consider $|I|$ sets of decision variables $z_i$. 

The last step is to take the (KKT) optimality conditions of the middle-level MPPDC problem \eqref{eq:middlelevel-mppdc} and add them to the top-level problem, resulting in a (trilevel) mathematical program with complementarity from primal and dual constraints. Similarly to the LCP-based reformulation summarised in Section \ref{sec:posit-semid-m-x-dep}, this strong duality reformulation has the requirement that the bottom level is not directly influenced by the middle-level decision variables. With a weak trilevel structure (as per Definition \ref{def:structure}), both the objective function and constraints are convex (or affine) and the KKT conditions of \eqref{eq:middlelevel-mppdc} are thus sufficient for optimality. However, to the best of our knowledge, no constraint qualification is known to hold for the problem \eqref{eq:middlelevel-mppdc}. For example, Slater's constraint qualification \parencite[all nonlinear constraints can be satisfied as strict inequalities,][]{slater1950lagrange} is not satisfied because weak duality states that 
$$
\sum_{i \in I} (z_i^\top F_i z_i + e_{i3}^\top z_i -  a_{i3}^\top p_i) \ge 0
$$ 
and thus, the nonlinear strong duality constraint \eqref{eq:strong-duality-mppdc} cannot be strictly satisfied. This means that the KKT conditions of this problem are only sufficient but not necessary for optimality. Nevertheless, this tells us that if we find a point that satisfies the KKT conditions, that point is optimal for the problem \eqref{eq:middlelevel-mppdc}. The complete single-level strong duality reformulation is thus
\begin{subequations}
  \label{eq:mpc-pdc}
  \begin{align}
    min_{x,y,z} \quad & c_1^\top x + d_1^\top y + e_1^\top z \\
    \st \quad & A_1 x + B_1 y + C_1 z \geq a_1 \\
    & 0 \le y \perp d_2 + B_2^\top \gamma \ge 0 \\
    & 0 \le \gamma \perp a_2 - A_2 x - B_{2} y - \sum_{i \in I} C_{i2} z_i \ge 0 \\
    & 0 \le z_i \perp e_{i2} + C_{i2}^\top \gamma - C_{i3}^\top p_i^b - F_i^\top z_i^b + (F_i+F_i^\top)z_i\epsilon + e_{i3}(x)^\top \epsilon \ge 0 \ \forall i \in I \label{eq:mpc-pdc-con4}\\
    & 0 \le p_i \perp C_{i3}z_i^b - a_{i3}(x)\epsilon \ge 0 \ \forall i \in I \label{eq:mpc-pdc-con5}\\
    & 0 \le z_i^b \perp e_{i3}(x) - C_{i3}^\top p_i + F_i z_i \ge 0 \ \forall i \in I \\
    & 0 \le p_i^b \perp C_{i3} z_i - a_{i3}(x) \ge 0 \ \forall i \in I \\
    & 0 \le \epsilon \perp - \sum_{i \in I} (z_i^\top F_i z_i + e_{i3}(x)^\top z_i -  a_{i3}(x)^\top p_i) \ge 0 \label{eq:mpc-pdc-strongduality},
  \end{align}
\end{subequations}
where $p_i^b$ and $z_i^b$ are the dual variables of the bottom-level primal and dual constraints \eqref{eq:middlelevel-mppdc-con2} and \eqref{eq:middlelevel-mppdc-con3}, respectively, and $\epsilon$ is the dual variable of the strong duality constraint \eqref{eq:strong-duality-mppdc}. That is, $p_i^b$ can be interpreted as the middle-level shadow prices associated with the bottom-level primal constraints. On the other hand, it is well known that the dual of the dual problem is the primal problem, and the dual variables associated with dual constraints are the primal variables. The value of these bottom-level primal variables must be the same for the middle- and bottom-level players, i.e., $z_i^b = z_i$. Note that the right-hand sides of constraints \eqref{eq:mpc-pdc-con4}, \eqref{eq:mpc-pdc-con5} and \eqref{eq:mpc-pdc-strongduality} contain bilinear terms (assuming $e_{i3}(x)$ and $a_{i3}(x)$ are affine) including the top-level variables $x$, making the resulting model nonconvex in general. As discussed before, such constraints can be modelled as quadratic SOS1 constraints and solved using spatial branch-and-bound-based methods.

If the problem instead has a strong trilevel structure, some of the terms $a_{i3}$ or $e_{i3}$ would effectively be functions of $y$, and the strong duality constraint \eqref{eq:strong-duality-mppdc} would consequently have nonconvex bilinear terms. The middle-level variables would be considered fixed for the bottom-level problems, but not for the middle level. A nonconvex strong duality constraint in the middle-level problem \eqref{eq:middlelevel-mppdc} would result in the KKT conditions of the problem not even being sufficient for optimality. If we assume for example that the middle-level variables $y$ appeared in linear terms added to the constant terms $a$ and $e$, constraint \eqref{eq:strong-duality-mppdc} would become 
$$ 
\sum_{i \in I}(z_i^\top F_i z_i + \left(e_{i3} + N_y^{obj} y \right)^\top z_i - \left(a_{i3}(x)+ N_y^{con} y \right)^\top p_i) \le 0, 
$$
resulting in bilinear terms $(N_y^{obj} y)^\top z$ and $(N_y^{con} y)^\top p$, where $N_y^{obj}$ and $N_y^{con}$ are the coefficient matrices of the $y$-variables in the bottom-level objective and constraints, respectively.

\subsection{Comparison of trilevel formulations}
\label{sec:comparison}

In the reformulation \eqref{eq:large-single-level-lpcc-in-x-and-y} \parencite{Gabriel_et_al:2021}, the vector $\tilde{z}$ contains both the primal and dual variables of each bottom-level problem. This is because the variable $\tilde{z}$ appears in the LCP \eqref{eq:lower-level-lcp}, which, in the problems presented in this paper, represents the concatenated KKT conditions of the bottom-level problems. We denote by $n_2$ the number of variables in the middle-level problem and by $m_2$ the number of constraints in the same problem, and analogously, $n_3$ and $m_3$ for the number of variables and constraints, respectively, in the bottom-level problem. There are then $n_2+m_2 + 4(n_3+m_3)$ complementarity constraints \eqref{subeq:large-single-level-lpcc-in-x-and-y.c}-\eqref{subeq:large-single-level-lpcc-in-x-and-y.h} and \eqref{subeq:large-single-level-lpcc-in-x-and-y.j} in formulation \eqref{eq:large-single-level-lpcc-in-x-and-y}, and one variable for each complementarity constraint. Additionally, there are $n_3+m_3+1$ equality constraints \eqref{subeq:large-single-level-lpcc-in-x-and-y.i} used in the reformulation, and the constraints \eqref{subeq:large-single-level-lpcc-in-x-and-y.b} for the top-level problem.

The novel strong duality formulation \eqref{eq:mpc-pdc} (assuming only inequality constraints and nonnegative variables in the middle- and bottom-level problems for comparison) results in $n_2+m_2$ complementarity constraints for the middle-level variables and constraints, $2(n_3+m_3)$ complementarity constraints for the bottom-level primal and dual variables and constraints, and one complementarity constraint for the strong duality. Because strong duality is represented as an inequality constraint, no equality constraints are needed for the strong duality reformulation. 

The strong duality reformulation of the bottom level results in half the number of complementarity constraints compared to the LCP reformulation presented in \textcite{Gabriel_et_al:2021}, plus one for strong duality, and no equality constraints. While the LCP reformulation results in two nested optimisation problems, the intermediate MPPDC \eqref{eq:middlelevel-mppdc} in the strong duality reformulation is a single problem, explaining the difference in the number of constraints. This is computationally beneficial, as large numbers of complementarity constraints contribute greatly to the computational challenges with equilibrium problems. Additionally, it should be noted that the column-and-constraint generation algorithm \parencite{dvorkin2017co} requires the middle- and bottom-level problems to be represented as a single optimisation problem, suggesting that the strong duality approach could be easily extended to that context, unlike the LCP solution set reformulation.

On the other hand, the main disadvantage of our strong duality formulation is that the strong duality constraint \eqref{eq:strong-duality-mppdc} retains the quadratic term from the bottom-level objective function, while the previous formulation has only affine constraints. This results in the formulation \eqref{eq:middlelevel-mppdc} not satisfying a constraint qualification, making the KKT conditions only sufficient for optimality. Additionally, unlike the strong duality formulation, the formulation in \textcite{Gabriel_et_al:2021} is applicable to settings where the bottom-level complementarity conditions are not derived as KKT conditions of an optimisation problem. For example, the spatial price equilibrium problem in \textcite{Gabriel_et_al:2021} could not be reformulated using strong duality. 

\section{Applications in energy-environmental planning}
\label{sec:applications}

In this section, we describe a trilevel power market equilibrium model that contains environmental considerations for the top-level regional policy-maker. Finding effective instruments for emission reduction and climate change mitigation is becoming increasingly important, and we focus our attention on carbon tax \parencite[see][for a review]{koppl2023carbon}. At the middle level, a single regional system operator is responsible for operating transmission lines $a \in A$ between nodes $k \in K$, maximising its profit from operating the system. 

At the bottom level, each energy producer $i \in I$ produces electricity at nodes $k \in K$ using energy sources $j \in J$ and sells the electricity to nodes $k' \in K$, that is, the electricity is not necessarily sold to the same node it is produced in. The producers maximise their profit from selling electricity, knowing that their decisions will affect the selling prices, making the bottom level a Cournot oligopoly. Instead of considering a fixed demand that must be satisfied exactly, we model the demand side as reacting with an affine relationship between production and price so that total demand increases linearly as the price of electricity decreases. This means, e.g., that if the producers started to generate unreasonably large amounts of electricity, the price would go down because more and more of the (elastic) demand is satisfied.

Finally, we consider a set $D$ of representative days \parencite{poncelet2016selecting} of renewable generation availability factors and demand curves. The top-level regulator chooses a tax and minimum renewable share which apply for all days. In contrast, the operational decisions at the system operator- and producer levels can differ between the days $d \in D$. The weights of the representative days are denoted with $P_d$, with $\sum_{d \in D}P_d = 1$, that is, $P_d$ represents the fraction of days in a year that is represented by day $d$. The purpose of representative days is to reduce the size and complexity of the model while still being able to realistically convey the variability in renewable energy availability and demand, and they are used in models such as US-REGEN \parencite{young2020regen} and LIMES-EU \parencite{nahmmacher2014limes}.

Our illustrative example is based on the model in \textcite{hobbs2001linear}. This model is chosen because of its simple nature, as using a more realistic model would require further discussion on assumptions and data, shifting the focus away from the methodological contributions of this paper. We highlight however, that While this reference model is a simplified representation of reality, models containing an equivalent structure to that in \textcite{hobbs2001linear} are used in case studies by, e.g., \textcite{keles2020cross, ribo2019effects}. 

\subsection{The top-level regulator}

On top of this trilevel hierarchy is the regional regulator which tries to both maximise the amount of electricity produced and minimise the carbon dioxide emissions from doing so. The motivation for this setting is to balance the utility from electricity generation and to maintain reasonable electricity prices, while simultaneously mitigating negative environmental outcomes. 

In addition to maximising production, the regulator wants to minimise the total emissions $\sum_{ijk} \eta_{ijk} z_{ijkd}$ from electricity generation, where $\eta_{ijk}$ is the emissions factor corresponding to the production level $z_{ijkd}$. For carbon-emitting energy sources, $\eta_{ijk} > 0$, while it is zero for zero-emission energy sources. These two objectives are then converted into a single objective by giving the total production value a weight $r \in (0,1)$ and the total emissions a weight $(1-r)$. By varying the value of this weight parameter, one could, for example, consider different priorities between these two objectives.

The top-level player decides on a carbon tax~$x$, which affects each firms' variable costs: $\gamma_{ijk} = \nu_{ijk} + \eta_{ijk}x$, where $\nu_{ijk} > 0$ is the cost specific to the firm-fuel combination $(i,j)$ in node $k$, and $\eta_{ijk}$ is the emissions factor.
Additionally, the top-level player can impose a minimum renewable share $\rho$ that the system operator must satisfy at each node $k \in K$. We assume $\rho$ to be the same for all nodes, but it would be straightforward to extend our model to consider this minimum renewable share to differ by node. The carbon tax and minimum renewable share affect the optimal solutions of the middle- and bottom-level players, resulting in different values for $z$, and consequently, the top-level objective value. Increasing the carbon tax results in lower emissions as the high-emission sources become more expensive for the producers. However, this also results in the market equilibrium in the lower levels shifting towards lower total production and higher electricity prices. 

Given the upper-level variables $x$ and $\rho$, the overall problem for this top-level player is given as
\begin{subequations}
\label{eq:regulator}
\begin{align}
   \max_{x,\rho,y,z} \quad &\sum_{d \in D} P_d \sum_{i \in I, j \in J, k \in K} (r - (1-r)\eta_{ijk}) z_{ijkd} \label{eq:regulator-obj}\\
  \st \quad &x, \rho \geq 0 \\
  &z \text{ and } y \text{ solve } \eqref{eq:profit_iso} \text{ for all } d \in D.
   \end{align}
\end{subequations}

\subsection{Profit-maximising system operator}
\label{sec:operator-problem}

At the middle level, following the model in \textcite{hobbs2001linear}, we consider a profit-maximising independent system operator (ISO). This ISO is responsible for operating the transmission lines $a \in A$ between nodes $k \in K$ for each representative day $d \in D$ and has to make sure that the lines function within their capacity limits, between $-T_a^-$ and $T_a^+$. The ISO chooses each node's net import $y_{kd}$ of electricity through the transmission lines (i.e., negative $y_{kd}$ implies that more electricity is produced than used in node $k$, and electricity is exported to other nodes). The line flows are determined from these using power transmission distribution factors (PTDFs) (see, e.g., \textcite{burrmetzler2000complementarity} for a thorough description). 

The ISO's problem for the representative day $d \in D$ can be stated as the following linear program.
\begin{subequations}
\label{eq:profit_iso}
\begin{align}
  \max_{y_{kd}, z_{ijkd}} \quad &\sum_{k \in K} w_{kd} y_{kd} \label{eq:profit_iso-obj}\\
  \st \quad &-\sum_{k \in K}PTDF_{ka} y_{kd} \le T_a^- \quad (\phi_{ad}^-)\ &&\forall a \in A \\
  &\sum_{k \in K}PTDF_{ka} y_{kd} \le T_a^+ \quad (\phi_{ad}^+)\ &&\forall a \in A \\
  &\sum_{i \in I, j \in R} z_{ijkd} \ge \rho\sum_{i \in I, j \in J} z_{ijkd} \quad (\psi_{kd})\ &&\forall k \in K \label{eq:profit_iso-rho} \\
  &z_{ijkd} \text{ solve } \eqref{eq:profit_max} \text{ for all } i \in I,
\end{align}
\end{subequations}
where $w_{kd}$ is a congestion-based wheeling fee for node $k \in K$ in day $d \in D$ and $R \subseteq J$ is the set of renewable energy sources. The wheeling fee is the unit price the producers have to pay to the ISO for selling electricity at node $k$, and the price that the ISO pays to the producer for each unit of electricity produced at node $k$, and the prices of buying and selling electricity in a node are assumed to be the same. The variables in parentheses to the right of each constraint are the corresponding dual variables.

Constraint \eqref{eq:profit_iso-rho} states that the ISO has to choose such transmission values that the renewable production share in each node is at least $\rho$, decided by the top-level regulator. We assume that the ISO has no mechanism for influencing the producers to, for example, increase their renewable share. This assumption results in a weak trilevel structure (Definition \ref{def:structure}). Instead of directly influencing the producers, the optimistic bilevel assumption described earlier results in the ISO ``choosing'' the best (in terms of \eqref{eq:profit_iso-obj}) equilibrium solution for the bottom-level problems that satisfies \eqref{eq:profit_iso-rho}.

\subsection{Oligopoly of the producers}
\label{sec:bilevel-energy-oligopoly}

We next consider the lower-level optimisation problems for a set of energy firms $i \in I = \set{1, \dotsc, n_F}$. We start by presenting these problems formulated for a bilateral market where electricity producers sell directly to consumers, which turns out to be the simpler case, and then proceed to add arbitragers to arrive at a POOLCO market model where the producers instead sell their electricity to a central auction. The POOLCO model more accurately represents the Nordic system and is thus used in the case study in Section \ref{sec:casestudy}. For a detailed discussion on different market types, we refer the reader to \textcite{ilic1998power}.

Let us first assume that at this lower level, these $n_F$ firms constitute the entire market. Each firm has a production capacity in some of the nodes $k \in K$ and can bilaterally sell their electricity directly to any of the nodes. For production, the producers have a set of energy sources $j \in J$. Our formulation for this producer level follows the ideas in \textcite{hobbs2001linear}.

In this first model without arbitragers, every firm $i \in I$ decides on its sales and production for each node $k \in K$ and day $d \in D$, taking into account linear inverse demand functions $p_{kd}(s_{1kd}, \dotsc, s_{n_Fkd}) = \alpha_{kd} - \beta_{kd}\sum_{i=1}^{N_F} s_{ikd}$ with price intercept $\alpha_{kd} > 0$ and slope $\beta_{kd} > 0$. These parameters are assumed to vary per day, representing the changes in demand. Recall that $s_{ikd}$ is the amount of electricity sold by producer $i$ to node $k$ in day $d$, and the market price at node $k$ thus depends on the sum of the sales of all firms into node $k$. 

Additionally, each producer $i \in I$ has maximum production levels~$z_{ijkd}^{\max}$ determined by their installed production capacity. For wind and solar power, the maximum production level depends on the representative day $d$. Each producing firm solves the profit-maximisation problem
\begin{subequations}
\label{eq:profit_max}
\begin{align}
  \max_{s_{ikd}, z_{ijkd}} \quad &\sum_{k \in K} \left( \left(\alpha_{kd} - \beta_{kd}\sum_{i' \in I} s_{i'kd}\right) s_{ikd} - \sum_{j \in J} \gamma_{ijk} z_{ijkd} - (s_{ikd}-z_{ijkd})w_{kd} \right) \label{eq:prod_profit}\\
  \st \quad &z_{ijkd} \le z_{ijkd}^{max} \quad (\lambda_{ijkd}) \ \forall j \in J, k \in K \label{eq:prod_capacity}\\ 
  &\sum_{k \in K}s_{ikd} = \sum_{j \in J, k \in K}z_{ijkd} \quad (\theta_{id}) \label{eq:prod_balance} \\
  &z_{ijkd}, s_{ikd} \geq 0,
\end{align}
\end{subequations}
where $\gamma_{ijk}$ is the marginal production cost for firm $i$ in node $k$ with fuel type $j$, composed as the sum of a firm-specific cost~$\nu_{ijk}$ and an emissions cost~$\eta_{ijk}x$, depending on the carbon tax $x$ determined by the regulator. 

The first term in \eqref{eq:prod_profit}, involving the sales variables $s_{ikd}$ represents the revenue from selling electricity to different nodes $k \in K$. The nodal price is $p_{kd} = \alpha_{kd} - \beta_{kd} \sum_{i \in I} s_{ikd}$. The cost of producing energy is $\gamma_{ijk}$. The producers pay a wheeling fee $w_{kd}$, which is determined by the transmission network congestion and paid to the ISO. In this hub-network model, the wheeling fee is also what the ISO pays the producers for producing extra energy in each node $k$.

Constraint \eqref{eq:prod_capacity} states that production cannot exceed capacity $z^{max}_{ijkd}$ and constraint \eqref{eq:prod_balance} states that for each producer, the total sales must equal total production. It is easy to see that the objective function \eqref{eq:prod_profit} is concave for $\beta_{kd} >0$ and the constraints are affine. Thus, the bottom-level problem \eqref{eq:profit_max} has the same structure as the quadratic problems discussed in Section \ref{sec:mpc-pdc}. 

Finally, we include a market-clearing constraint
\begin{equation}
    \sum_{i \in I} s_{ikd} - \sum_{i \in I, j \in J} z_{ijkd} = y_{kd} \quad (w_{kd})\ \forall k \in K, d \in D. \label{eq:market_clearing}
\end{equation} 
This constraint is similar to constraint \eqref{eq:prod_balance}, which instead considers the difference between sales and production for each producer $i \in I$. 
We adopt the Bertrand assumption used in \textcite{hobbs2001linear}: the system operator sees the wheeling fees as fixed, instead of using market power to affect their values. In order to achieve this, the market-clearing constraint \eqref{eq:market_clearing} is considered outside the system operator and producer problems, appearing ``separately'' in the final single-level formulation, effectively becoming a top-level constraint.

\subsubsection{Extending the producer oligopoly: including arbitrage}
\label{sec:arb}

We are interested in modelling the Nordic market and, to achieve that, we extend the bilateral market model represented by \eqref{eq:profit_max} into a POOLCO model. In a POOLCO market model, it is assumed that the producers sell their electricity to a central auction where the price is determined based on the amount of sold electricity and network congestion. \textcite{burrmetzler2000complementarity} and \textcite{hobbs2001linear} show that a bilateral market with \emph{arbitragers} is equivalent to a POOLCO market, assuming Cournot competition. Arbitragers are bottom-level players who have no production capacity, but they instead make their profits by exploiting the price differences between nodes, buying cheap electricity and selling it to nodes with a higher price. They act as price-takers and thus do not anticipate their effect on the price $p_{kd}$. The arbitrager's problem is
\begin{subequations}
\label{eq:arb}
\begin{align}
  \max_{a} \quad &\sum_{k \in K} (p_{kd} - w_{kd}) a_{kd} \\
  \st \quad &\sum_{k \in K} a_{kd} = 0 \quad (p_d^H),
\end{align}
\end{subequations}
where $a_{kd}$ is the amount of electricity sold by the arbitrager to node $k$ in day $d$ and the price at node $k \in K$ depends on the sales from the producers and the arbitragers, thus becoming $p_{kd} = \alpha_{kd} - \beta_{kd}\left(\sum_{i \in I, j \in J} s_{ikd} + a_{kd}\right)$. We can trivially obtain the KKT conditions of \eqref{eq:arb}, a linear maximisation problem (recall that the arbitragers are price-takers, and $p_{kd}$ is thus treated as a constant). The KKT conditions \eqref{eq:bottom-level-primal-arbitrage.con3} and \eqref{eq:bottom-level-primal-arbitrage.con4} are necessary and sufficient for optimality and adding them to \eqref{eq:profit_max}, we obtain
\begin{subequations}
\label{eq:bottom-level-primal-arbitrage}
\begin{align}
  \max_{s_{ikd}, z_{ijkd}, a_{ikd}} \quad &\sum_{k \in K} \left(\left(\alpha_{kd} - \beta_{kd} \left(\sum_{i' \in I}s_{i'kd} + a_{ikd}\right) \right) s_{ikd} - \sum_{j \in J} \gamma_{ijk} z_{ijkd}  - (s_{ikd}-z_{ijkd})w_{kd} \right) \label{eq:bottom-level-primal-arbitrage.obj}\\
  \st \quad &z_{ijkd} \le z_{ijkd}^{max} \quad (\lambda_{ijkd}) \ \forall j \in J, k \in K \label{eq:bottom-level-primal-arbitrage.con1}\\ 
  & \sum_{k \in K}s_{ikd} = \sum_{j \in J, k \in K}z_{ijkd} \quad (\theta_{id}) \label{eq:bottom-level-primal-arbitrage.con2}\\
  & \alpha_{kd} - \beta_{kd}\left(\sum_{i' \in I}s_{i'kd} + a_{ikd}\right) =p^H_{id} + w_{kd} \ \forall k \in K \label{eq:bottom-level-primal-arbitrage.con3}\\
  & \sum_{k \in K} a_{ikd} = 0 \label{eq:bottom-level-primal-arbitrage.con4}\\
  & z_{ijkd}, s_{ikd} \geq 0, \label{eq:bottom-level-primal-arbitrage.vars}
\end{align}
\end{subequations}
where $a_{ikd}$ is the net amount of power sold in node $k$ by the arbitrager(s), and $p^H_{id}$, the dual variable associated with the arbitrager constraint, is the price at the central auction $H$. Both $a_{ikd}$ and $p^H_{id}$ are indexed over the different producers $i \in I$, to highlight that each producer can influence these values with their decisions, and to avoid decision variables shared by players. This would result in a generalised Nash equilibrium problem \parencite{facchinei2010generalized} that would be computationally more challenging. However, the values $a_{ikd}$ and $p^H_{id}$ are the same for all producers at equilibrium, as shown in Appendix \ref{app:substitutions}, and the approach of having separate variables for each producer is thus valid. Constraint \eqref{eq:bottom-level-primal-arbitrage.con3} can be therefore written as $p_{kd} - w_{kd} = p^H_{id}$. That is, including arbitragers results in the producers selling their electricity to the central auction at the hub price $p^H_{id}$ (or simply $p^H_d$ at equilibrium), which is the sum of the price $p_{kd}$ at node $k \in K$ and the wheeling fee $w_{kd}$ paid to the system operator. Constraint \eqref{eq:bottom-level-primal-arbitrage.con4} states that since the arbitragers have no production capacity, their net sales amounts must be zero. The objective function is still concave after adding the arbitrage variables, and the new constraints are affine. \textcite{burrmetzler2000complementarity} shows further substitutions and simplifications to the producer and system operator problems, which are shown in Appendix \ref{app:substitutions}, along with the resulting model that is used for the computational experiments in Section \ref{sec:comp}.

\section{Computational experiments}
\label{sec:comp}

To illustrate the performance of the trilevel optimisation framework in a realistic problem setting, we solve the trilevel model described in the previous section, using randomly generated instances of varying sizes. The data used in these computational experiments mimics the data in the case study of \textcite{belyak2023optimal}, whose data is from the ENTSO-E Transparency Platform \parencite{hirth2018entso} and is further described in Section \ref{sec:casestudy}. 
The computational experiments were performed using 8 CPU threads and 16GB of RAM. All code was implemented in Julia v1.7.3 \parencite{Julia-2017} using the Gurobi solver v10.0.0 \parencite{gurobi} and JuMP v1.5.0 \parencite{DunningHuchetteLubin2017} and is available in \textcite{repo}.

\subsection{Comparing formulations}
\label{sec:formulations-comp}

We compare the performance of the two single-level reformulations, the LCP-based reformulation from \textcite{Gabriel_et_al:2021} (Section \ref{sec:posit-semid-m-x-dep}) and our strong duality reformulation (Section \ref{sec:mpc-pdc}) by solving 50 randomly generated problems with 2 producers, 5 energy sources, 3 nodes and 3 representative days. This problem size was chosen as the base case because it seems to be large enough to make the problems challenging to solve, but small enough for them to be mostly solvable within a time limit of one hour.

\begin{figure}[!ht]
    \centering
    \includegraphics[width=0.85\textwidth]{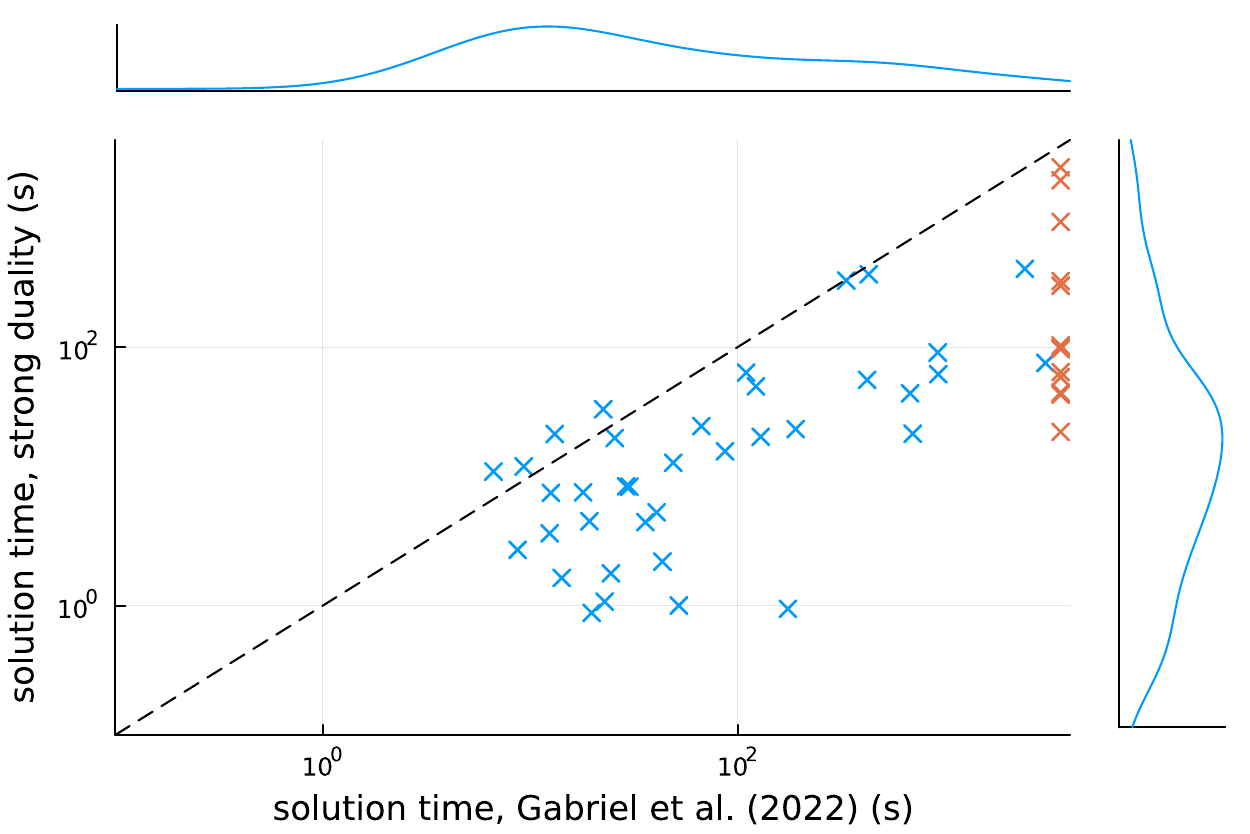}
    \caption{Solution times for the two formulations on 50 random instances with 2 producers, 5 energy sources, 3 nodes and 3 representative days. If one of the methods failed to find a solution within 3600s, a red marker is used, and the marginal distributions on the right and top sides exclude unsolved instances.}
    \label{fig:soltime_scatter}
\end{figure}

The results are presented in Figure \ref{fig:soltime_scatter} and the main observation here is that the novel strong duality formulation is faster in most cases. In Figure \ref{fig:soltime_scatter}, markers below the diagonal (dashed line) correspond to such cases. In 13 instances, the formulation of \textcite{Gabriel_et_al:2021} did not find an optimal solution in an hour while our strong duality formulation did. One major issue with both models compared here is that usually the first feasible solutions are found at the end of the solution process and most of the solution time is spent on improving the dual bound without finding any feasible solutions. Nevertheless, changing solver parameters to emphasize finding feasible solutions was not found to have a major impact on performance.

As discussed in Section \ref{sec:comparison}, the strong duality formulation results in fewer constraints than the reformulation in \textcite{Gabriel_et_al:2021}. Recall that in our models, complementarity constraints are formulated as SOS1 constraints. The model sizes in the base case test problems are presented in Table \ref{tbl:modelsizes}, and out of the two, our strong duality model is smaller, except for having one more quadratic SOS1 constraint to represent the strong duality constraint \eqref{eq:mpc-pdc-strongduality}. In our model, all non-complementarity quadratic constraints are inequality constraints, while in the LCP-based model, there is one quadratic equality constraint (the first one in \eqref{subeq:large-single-level-lpcc-in-x-and-y.i}).

\begin{table}[!ht]
\centering
\begin{tabular}{r|cc}
            & strong duality (this paper) & LCP \parencite{Gabriel_et_al:2021} \\ \hline 
variables   & 678     & 949  \\
affine constraints  & 306     & 757  \\
quadratic constraints & 100     & 100  \\
affine SOS1   & 288     & 648  \\
quadratic SOS1  & 100     & 99  
\end{tabular}
\caption{Model sizes for the two reformulations}
\label{tbl:modelsizes}
\end{table}

\begin{figure}[!ht]
    \centering
    \includegraphics[width=0.85\textwidth]{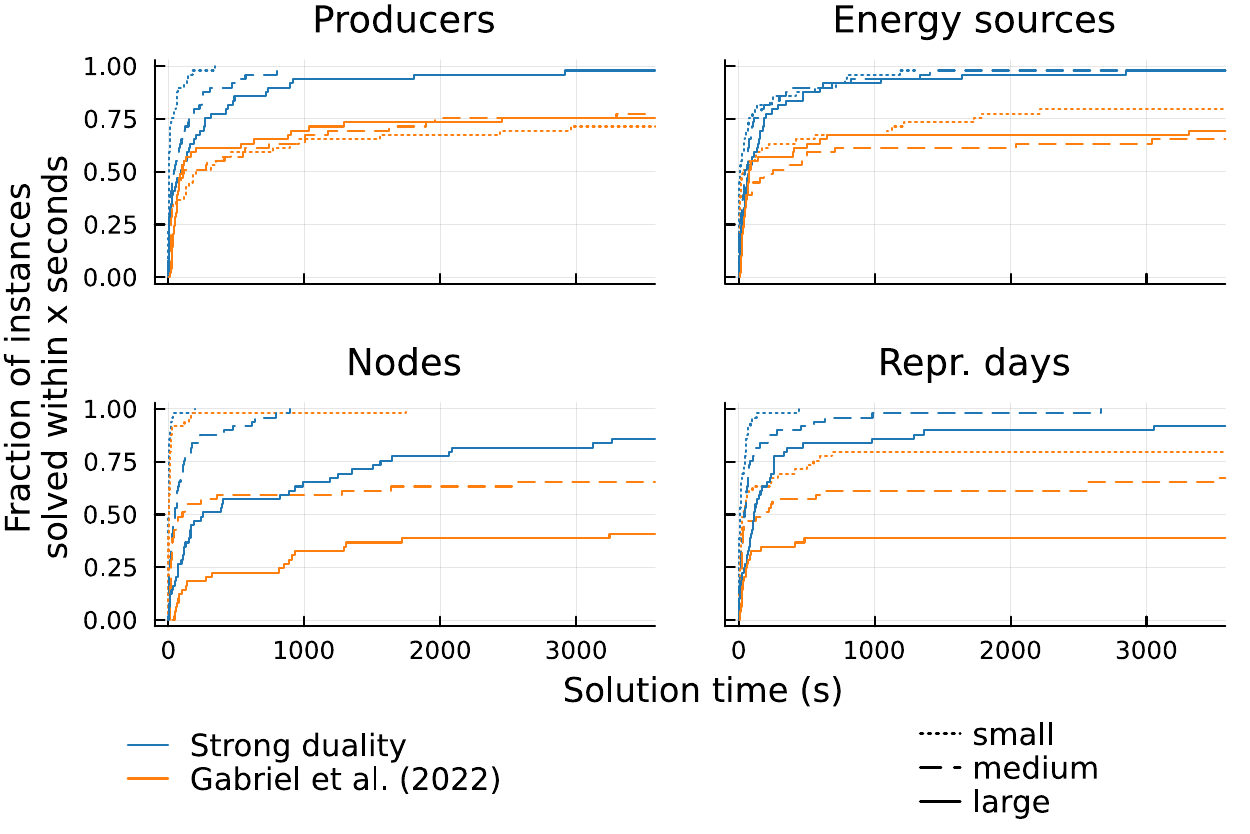}
    \caption{Cumulative distribution functions of solution times for the two formulations with 1-3 producers, 4-6 energy sources and 2-4 nodes and representative days. For each problem size, 50 instances are generated and solved.}
    \label{fig:performanceprofiles}
\end{figure}

Next, we analyze how problem size affects solution times by varying either the number of producers, energy sources, nodes and representative days from the base case, one parameter at a time. The results are presented in Figure \ref{fig:performanceprofiles}. The medium-sized cases in each subfigure are similar to each other, which is expected as the problem sizes are the same. Varying the number of producers or energy sources seems to have only a small effect on the solution times while changing the number of nodes has a far stronger effect. The effect of the number of representative days is stronger than that of the number of producers and energy sources but seems to be weaker than that of the number of nodes. 

We can also see that the number of problems that were not solved to optimality within the time limit is affected by the number of nodes and representative days, but not by the number of firms or energy sources. Additionally, the novel strong duality model finds an optimal solution more frequently than the previous formulation. As predicted in Section \ref{sec:mpc-pdc}, the larger number of complementarity constraints in the LCP formulation (Table \ref{tbl:modelsizes}) proves to be computationally challenging, and the smaller strong duality model is solved faster.

\subsection{Case study: a five-node Nordic energy system}
\label{sec:casestudy}

The case study in \textcite{belyak2023optimal} considers five nodes, representing Finland, Sweden, Norway, Denmark and the combined Baltic countries (Estonia, Latvia and Lithuania). There are five producers, each owning production capacity in one of the five nodes. Nine different energy sources are available, consisting of five conventional sources: nuclear, coal, gas (closed- and open-cycle) and biomass, and four renewable sources: solar, hydro, onshore and offshore wind. Additionally, we consider three representative days of renewable generation availability factors and demand curves. Recall that in our model, the top-level regulator makes their decisions independent of the day considered, that is, the carbon tax and minimum renewable share are constant across different representative days. These representative days are obtained in \textcite{belyak2023optimal} by performing hierarchical clustering on demand, price and renewable availability data.

Day 1 is a winter day with higher demand, low solar availability and medium wind availability. Days 2 and 3 have a lower demand with day 2 representing a windy day with medium solar availability, and day 3 representing a sunny day with low wind availability. The details of the hierarchical clustering process can be found in \textcite{belyak2023optimal}.

\begin{figure}[!ht]
\centering
     \subfloat[][No carbon tax.]{
         \includegraphics[width=0.45\textwidth]{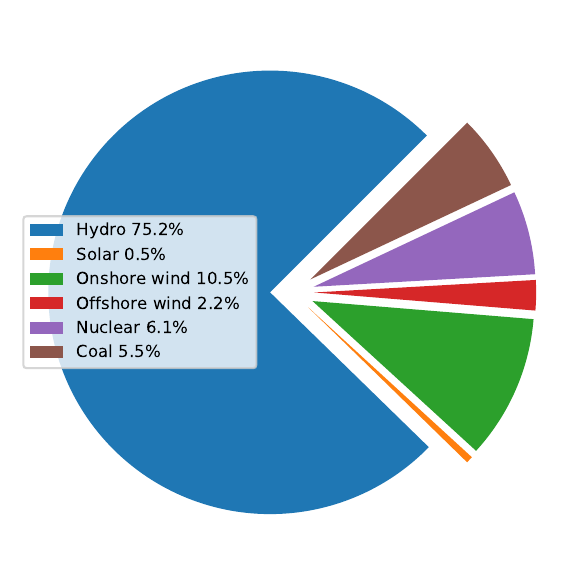}
         \label{fig:pie_notax}
     }%
     \qquad
     \subfloat[][Carbon tax 23e/ton.]{
         \includegraphics[width=0.45\textwidth]{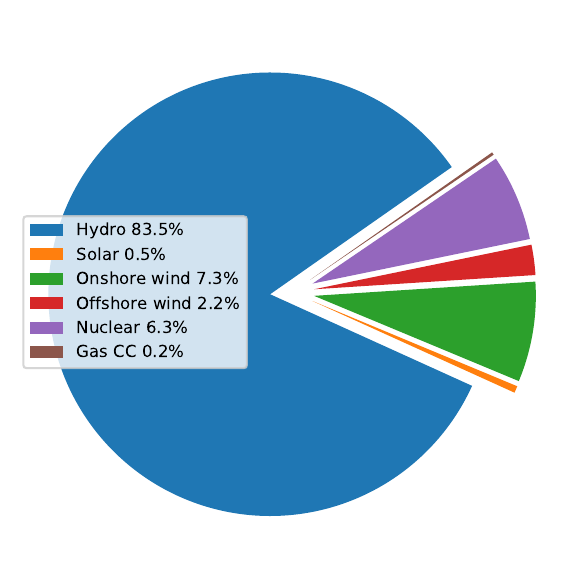}
         \label{fig:pie_tax}
     }
     \caption{Weighted average electricity production portfolio over the five nodes and three representative days.}
     \label{fig:pie}
\end{figure}

In Figure \ref{fig:pie}, the production portfolio (a weighted average over the representative days) is presented for a model with no carbon tax (i.e., the regulator heavily prefers maximising production over minimising emissions) and a carbon tax of 23 \euro/ton (enough to remove nearly all emissions). Compared to the baseline with no carbon tax, this 23 \euro/ton tax decreases the total production by 2.8\%. In this example, these carbon tax values are achieved by setting the weight parameter $r$ in the top-level objective \eqref{eq:regulator-obj} to 0.8 and 0.4, respectively. 

\begin{figure}[!ht]
    \centering
    \includegraphics[width=0.9\textwidth]{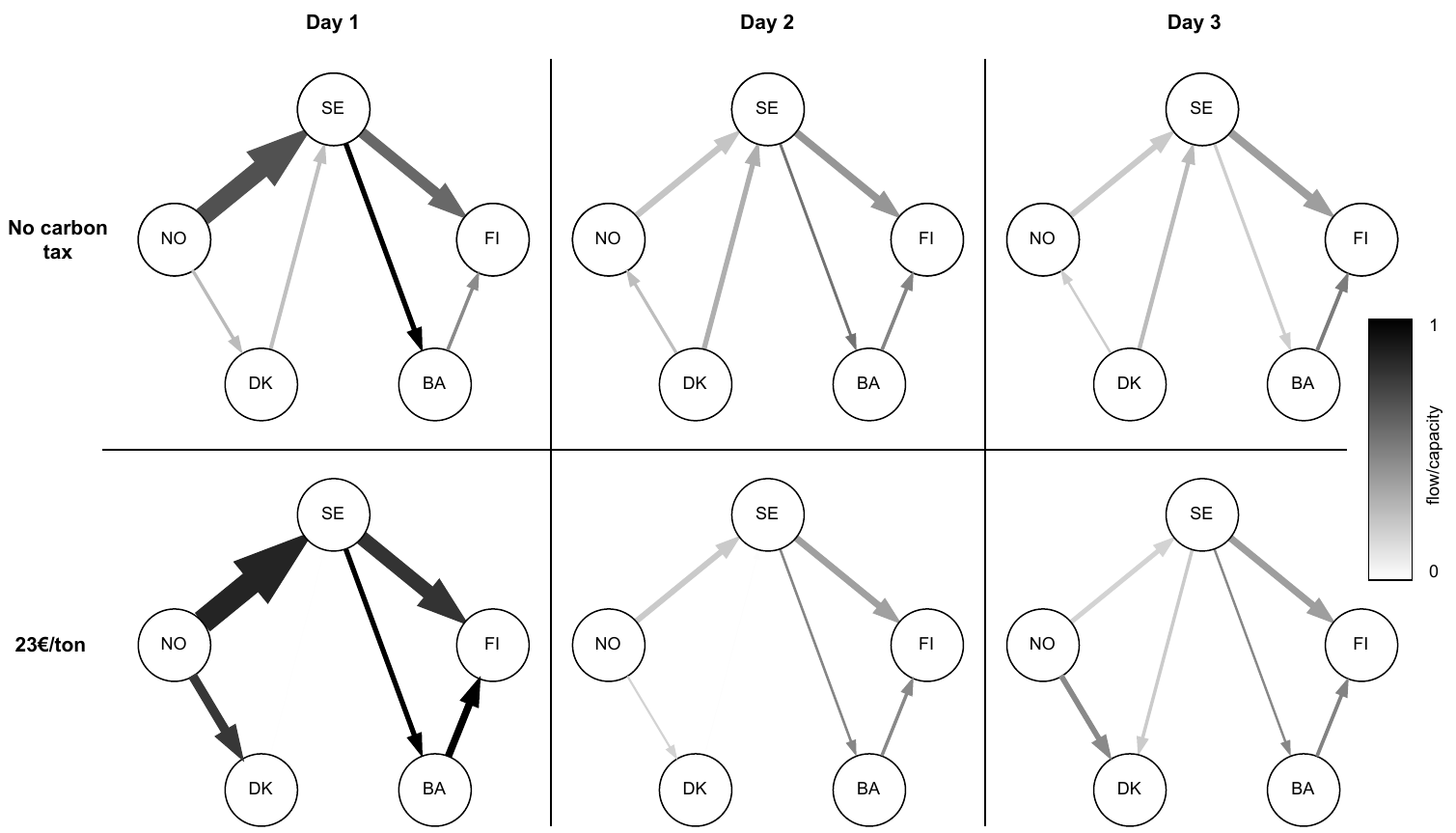}
    \caption{Transmission grid usage with different carbon taxes and representative days. The size of an arrow is proportional to the flow on the line and the color of an arrow represents congestion: black arrows correspond to lines operating at their limit. The nodes are FI=Finland, SE=Sweden, NO=Norway, DK=Denmark, BA=Baltic countries.} 
    \label{fig:flows}
\end{figure}

Because of the substantial hydropower production capacity in the Nordic system, particularly in Norway and Sweden \parencite{irena2023}, the renewable share of production is large even without a carbon tax. A part of the increase in hydropower usage when the carbon tax is introduced comes from decreasing onshore wind production. This is a consequence of the simplified nature of the model and data, as the operational costs for both hydropower and onshore wind power are zero, resulting in multiple optima and indifference for the producers to use one or the other, as long as the production capacity of neither is exceeded. This artefact of the model could be easily removed by, e.g., setting the operational cost of either energy source to a small positive number instead of zero, causing the producers to prefer the cheaper source. However, this would imply an artificial preference for one source over the other. The only significant source of emissions is coal, and introducing a carbon tax of 23 \euro/ton removes all coal from the portfolio, bringing in a small amount of closed-cycle gas power instead. The closed-cycle gas production occurs in the Baltics for day 1, and to understand this emergence of gas better, we must examine the transmission network in Figure \ref{fig:flows}.

The first representative day has the highest network usage with large amounts of electricity transmitted from Norway to Finland through Sweden. With the carbon tax, the importance of transmission is further highlighted as the hydropower capacity in Norway is used for lowering overall prices under high demand and low production from both solar power and high-emission sources. The differences between representative days 2 and 3 are more subtle, but we can see, e.g., the reliance on wind power in Denmark: in the low-wind day 3, the carbon tax results in Denmark importing a significant amount of electricity from Norway, compared to the high-wind day 2. On the first day with a carbon tax of 23 \euro/ton, both lines connecting the Baltic countries to the rest of the system are at their capacity, explaining why the Baltic countries start using gas power after a carbon tax is introduced. This illustrates the complex interplay between the three levels that is captured by our model.

\section{Conclusions}
\label{sec:conclusions}

In this paper, we propose a novel formulation for trilevel optimisation problems focusing on energy systems planning with environmental considerations. Additionally, we characterise the notion of weak and strong trilevel structures and compare the computational performance of the novel strong duality-based reformulation in this paper and the LCP-based reformulation in \textcite{Gabriel_et_al:2021}.

The computational results are encouraging, as we are able to solve the case study to optimality within a few minutes despite the fact that both single-level reformulations considered are nonconvex problems. However, we note that preliminary experiments with seemingly small extensions to the model, such as adding ramping constraints (limiting the change in production between consecutive periods) to the producer problem made the problem computationally intractable. The small size of the case study is indicative of the very challenging (nonconvex) nature of these problems, and the authors note that the reformulations and solution methods in this paper should be viewed as one of the first steps towards an efficient solution framework for tri-level problems.

For the results in this paper, an off-the-shelf solver is used, which is useful to ensure a low barrier-to-entry for using the developed formulation. However, we believe that the computational performance can be increased considerably using specialised solution methods like column-and-constraint generation \parencite{dvorkin2017co}. Notably, ideas such as bilevel branch-and-bound \parencite{fischetti2018use} and convex hull reformulations of the middle-level feasible region \parencite{santana2020convex} may be explored in the context of the problems presented in this paper. In addition, the model could also be extended to consider transmission and/or production capacity expansion over multiple time periods, especially if computationally more efficient reformulations and solution methods are developed.

Despite the outstanding computational challenges, we show that the novel reformulation improves computational performance compared to the previous formulation \parencite{Gabriel_et_al:2021}, and we show that the framework can be applied to a setting representing the Nordic electricity market, and results on the effect of carbon tax can be obtained. A limitation of the formulation approach presented in this paper and that originally proposed by \cite{Gabriel_et_al:2021} is that they require a weak trilevel structure. In practice, relevant problems may instead have a strong trilevel structure, precluding the use of these reformulations. Thus, further research is needed on developing (heuristic) solution methods for problems with a strong trilevel structure.

\section*{Acknowledgements}
The calculations presented in this paper were performed using computer resources within the Aalto University School of Science “Science-IT” project.

Fabricio Oliveira and Olli Herrala were supported by the Research Council of Finland (decision numbers 348094 and 332180). Steven A. Gabriel was supported by the National Science Foundation (NSF) Award \#2113891, Civil Infrastructure Systems. Tommi Ekholm was supported by the Research Council of Finland(decision number 341311).

\printbibliography

\appendix

\section{Reformulation of a bottom-level LCP with a positive semi-definite M}
\label{app:posit-semid-m-x-dep}

\textcite{Cottle_et_al:2009} shows that if the matrix $M$ is positive semidefinite, all solutions to the LCP

\begin{equation} 
0 \leq z \perp q + N_x x + N_y y + Mz \geq 0,
\end{equation}
can be obtained as the following polyhedral set:

\begin{align}
  \{z \in R^{n_z}_{\geq 0} : & \ q + N_x x + N_y y + Mz \geq 0, \nonumber\\
  & \ (q + N_x x + N_y y)^T (z - \bar{z}) = 0, \label{eq:lcp-solution-set} \\
  & \ (M + M^T)(z - \bar{z}) = 0\}, \nonumber
\end{align}
where $\bar{z}$ is a solution to the LCP.

Hence, the middle-level problem can be re-written as

\begin{subequations}
\label{eq:lower-level-lpcc-reform}
  \begin{align}
    min_{y,z \geq 0} \quad & d_2^\top y + e_2^\top z \label{eq:middle-level-obj}\\
    s.t. \quad & A_2 x + B_2 y + C_2 z \geq a_2 \\
    & q + N_x x + N_y y + Mz \geq 0 \label{eq:lpcc-sol_set-1} \\
    & (q + N_x x + N_y y)^\top (z - \bar{z} )= 0 \label{eq:lpcc-sol_set-2} \\
    & (M + M^\top) (z - \bar{z} ) = 0. \label{eq:lpcc-sol_set-3}
  \end{align}
\end{subequations}

We observe that \eqref{eq:lpcc-sol_set-2} includes a bilinear term $y^\top N_y^\top z$ in an equality constraint. This is a nonconvex constraint, precluding the direct use of KKT conditions for obtaining an optimal solution to \eqref{eq:lower-level-lpcc-reform}. However, for problems with a weak tri-level structure, $N_y = 0$ and these bilinear terms vanish. In the next theorem, we assume $N_y = 0$. 

\begin{theorem}
  Let $M$ be a positive semidefinite matrix. Then, $(x^*, y^*, z^*)$ is an optimal solution of Problem~\eqref{eq:toplevel-prob} with middle level~\eqref{eq:middlelevel-prob} if and only if $(x^*, y^*, z^*, \bar{z}^*)$ is an optimal solution of the problem
  
  \begin{subequations}%
    \label{eq:bilevel-lpcc-in-x-and-y}%
    \begin{align}
      min_{x, y, z, \bar{z}} \quad & c_1^\top x + d_1^\top y + e_1^\top z \label{subeq:bilevel-lpcc-in-x-and-y.a}\\
      s.t. \quad & A_1 x + B_1 y + C_1 z \geq a_1 \label{subeq:bilevel-lpcc-in-x-and-y.b} \\
      & \bar{z} \in \argmin_{z^\prime \geq 0} \{ {z^\prime}^\top (q + N_x x + M{z^\prime}): \nonumber \\
      & \qquad \qquad \quad q + N_x x + M {z^\prime}\geq 0 \quad (\beta)\}, \label{subeq:bilevel-lpcc-in-x-and-y.c-d-e} \\
      & y,z \in \argmin_{\hat{y},\hat{z} \geq 0} \{d_2^\top \hat{y} + e_2^\top \hat{z} : \nonumber \\
      & \qquad \qquad  q + N_x x + M\hat{z} \geq 0 \quad (\delta) \nonumber \\
      & \qquad \qquad  (q + N_x x)^\top (\hat{z} - \bar{z}) = 0 \quad (\zeta) \label{subeq:bilevel-lpcc-in-x-and-y.f-g-h-i} \\
      & \qquad \qquad  (M + M^\top) (\hat{z}-\bar{z})=0 \quad (\eta) \nonumber \\
      & \qquad \qquad  A_2 x + B_2 \hat{y} + C_2 \hat{z} \geq a_2 \quad (\gamma) \}, \nonumber\\
      &\text{such that } (\bar{z}^*)^\top (q + N_x x^* + M\bar{z}^*) = 0.\nonumber
    \end{align}
  \end{subequations}
\end{theorem}
See Theorem 6 in \textcite{Gabriel_et_al:2021} for a proof of this result as well as related theoretical aspects of the general form of the problem.

The two nested optimization problems in~\eqref{eq:bilevel-lpcc-in-x-and-y} are a convex QP \eqref{subeq:bilevel-lpcc-in-x-and-y.c-d-e} and an LP \eqref{subeq:bilevel-lpcc-in-x-and-y.f-g-h-i}. Hence, the KKT conditions of both problems are necessary and sufficient for optimality and the two inner problems can be replaced by their necessary and sufficient KKT conditions, leading to the single-level reformulation \eqref{eq:large-single-level-lpcc-in-x-and-y}.

\section{Formulating the dual of a QP with affine constraints}
\label{app:dual-formulation}

Given a quadratic program with affine constraints
\begin{subequations}
\begin{align}
  \min_{z_i} \quad & \frac{1}{2} z_i^\top F_i z_i + e_{i3}(x)^\top z_i \\
  \st \quad & C_{i3} z_i \ge a_{i3}(x) \quad (p_i) \\
  & z_i \geq 0 \quad (s_i),
\end{align}
\end{subequations}
where we assume $F_i$ is a positive semidefinite symmetric matrix, the Lagrangian of the problem is 

\begin{equation}
    L(z_i,s_i,p_i) = \frac{1}{2} z_i^\top F_i z_i + e_{i3}(x)^\top z_i + p_i^\top(a_{i3}(x) - C_{i3} z_i) - s_i^\top z_i, \label{eq:lagrangian-first}
\end{equation}
where $p_i$ and $s_i$ are nonnegative Lagrange multipliers or dual variables. The first-order optimality condition is thus

\begin{equation}
    \nabla_{z_i} L(z_i,s_i,p_i) = F_i z_i + e_{i3}(x) - C_{i3}^\top p_i - s_i = 0
\end{equation}
%
and rearranging \eqref{eq:lagrangian-first} gives us 

\begin{equation}
    \frac{1}{2} z_i^\top F_i z_i - z_i^\top C_{i3}^\top p_i+ z_i^\top e_{i3}(x) - z_i^\top s_i + a_{i3}(x)^\top p_i, \label{eq:lagrangian-intermediate}
\end{equation}
which, using the first order condition $-F_i z_i = -C_{i3}^\top p_i + e_{i3}(x) - s_i$, becomes



\begin{equation}
    -\frac{1}{2} z_i^\top F_i z_i + a_{i3}(x)^\top p_i \label{eq:lagrangian-final}.
\end{equation}
Maximising Eq. \eqref{eq:lagrangian-final}, subject to the first-order optimality condition for $z_i$ and treating $s_i$ as a slack variable and removing its explicit representation from the problem results in the Lagrangian dual formulation

\begin{subequations}
\begin{align}
  \max_{p_i, z_i} \quad & -\frac{1}{2} z_i^\top F_i z_i + a_{i3}(x)^\top p_i \\
  \st \quad & C_{i3}^\top p_i - F_i z_i \le e_{i3}(x)\\
  & p_i \geq 0.
\end{align}
\end{subequations}

\section{Further substitutions for the middle and bottom levels}
\label{app:substitutions}

The producer model can be further simplified using the substitution $s_{ikd} = \sum_{j \in J} z_{ijkd}$, removing the sales variables and the balance constraint \eqref{eq:bottom-level-primal-arbitrage.con2}. For a further reduction, the remaining equality constraints \eqref{eq:bottom-level-primal-arbitrage.con3} and \eqref{eq:bottom-level-primal-arbitrage.con4} can be used to solve for $a$ and $p^H$. 

We have the necessary and sufficient KKT conditions 
\begin{align}
    & \alpha_{kd} - \beta_{kd}(\sum_{i' \in I}s_{i'kd} + a_{ikd}) = p^H_{id} + w_{kd} \forall k \in K \\
    & \sum_{k \in K} a_{ikd} = 0\\
\end{align}
of the arbitrager's problem, and with the substitution $\sum_{j \in J} z_{ijkd} = s_{ikd}$, we get 

\begin{align}
    & \alpha_{kd} - \beta_{kd}(Z_{kd} + a_{ikd}) = p^H_{id} + w_{kd} \forall k \in K \\
    & \sum_{k \in K} a_{ikd} = 0, \\
\end{align}
where $Z_{kd} = \sum_{i \in F, j \in J} z_{ijkd}$. In matrix form, we get

\begin{equation}
    \begin{bmatrix}
        Q_d & \mathbf{1} \\
        \mathbf{1}^\top & 0
    \end{bmatrix}
    \begin{bmatrix}
        a_{id} \\ p^H_{id}
    \end{bmatrix}
    = 
    \begin{bmatrix}
        \alpha_d - QZ_d - w_d \\ 0
    \end{bmatrix},
\end{equation}
where $Q_d$ is a square diagonal matrix with the element on the $k$th row and column being $\beta_{kd}$ and $\mathbf{1}$ is a vector of ones. It can be shown that 

\begin{equation}
    \begin{bmatrix}
        Q_d & \mathbf{1} \\
        \mathbf{1}^\top & 0
    \end{bmatrix}^{-1}
    = 
    \begin{bmatrix}
        L^d & h_d \\
        h^\top_d & \hat{h}_d
    \end{bmatrix},
\end{equation}
where 

\begin{align*}
    \hat{h}_d &= \frac{1}{\sum_{k \in K} \beta_{kd}^{-1}} \\
    h_{kd} &= \beta_{kd}^{-1}\hat{h}_d \\
    L^d_{k,k} &= \hat{h}_d\beta_{kd}^{-1}\sum_{k' \in K \setminus k}\beta_{kd}'^{-1} \\
    L^d_{k,k'} &= -\hat{h}_d\beta_{kd}^{-1}\beta_{kd}'^{-1}, \ k \neq k'.
\end{align*}
This results in the solution

\begin{align}
    a_{ikd} &= h_{kd} Z_d - Z_{kd} - \sum_{k' \in K} (\alpha_{k'd}-w_{k'd})L^d_{k,k'} \\ 
    p^H_{id} &= \sum_{k \in K} (\alpha_{kd} - w_{kd})h_{kd} - Z_d \hat{h}_d,
\end{align}
where 

\begin{align*}
    Z_d &= \sum_{i \in I, j \in J, k \in K} z_{ijkd} \\
    Z_{kd} &= \sum_{i \in I, j \in J} z_{ijkd}.
\end{align*} 

It can be seen that the values of $a_{ikd}$ and $P^H_{id}$ are the same for each firm $i \in I$ and we can drop the index $i$. \textcite{burrmetzler2000complementarity} shows that the arbitrage amounts correspond to the transmission values: $a_{kd}=y_{kd}$. 

These substitutions result in the problem formulation
\begin{subequations}
\label{eq:bottom-level-primal-arbitrage-reduced}
\begin{align}
  \max_{z_{id}} \quad & \left( \sum_{k \in K} (\alpha_{kd} - w_{kd})h_{kd} - Z_d \hat{h}_d \right) Z_{id} - \sum_{j \in J, k \in K} \left(\gamma_{ijk}-w_{kd}\right) z_{ijkd} \label{eq:bottom-level-primal-arbitrage-reduced.obj}\\
  \st \quad &z_{ijkd} \le z_{ijkd}^{max} \quad (\lambda_{ijkd}) \ \forall j \in J, k \in K \label{eq:bottom-level-primal-arbitrage-reduced.con1}\\ 
  & z_{ijkd}\geq 0 \label{eq:bottom-level-primal-arbitrage-reduced.vars}
\end{align}
\end{subequations}
and we can see that the substitutions do not change the concavity of the objective function: the quadratic term for producer $i \in I$ is $\hat{h}_d Z_{id}^2$. Finally, the sales variables $s_{ikd}$ are also eliminated from the market-clearing constraint, resulting in
\begin{equation}
h_{kd} Z_d - Z_{kd} + \sum_{k' \in K} (\alpha_{k'd} - w_{k'd})L^d_{kk'} = y_{kd} \quad (w_{kd})\ \forall k \in K. \label{eq:market_clearing_arb}
\end{equation}

The formulation \eqref{eq:bottom-level-primal-arbitrage-reduced} can be converted into an LCP by using the KKT optimality conditions. The combined KKT conditions of \eqref{eq:bottom-level-primal-arbitrage-reduced} for all producers $i \in I$ are
\begin{subequations}
\label{eq:bottom-level-primal-KKT}
\begin{align}
  0 \le z_d \perp &B_d z_d + \lambda_d + q^z_d \ge 0\\ 
  0 \le \lambda_d \perp & - z_d + z^{max}_d \ge 0,
\end{align}
\end{subequations}
where $q^z_{ijkd} = -\sum_{k \in K}(\alpha_{kd}-w_{kd})h_{kd} + (\gamma_{ijk} - w_{kd})$ and $B_d$ is a positive semidefinite 
matrix with $$B_d(ijk, i'j'k') = \begin{cases}
2\hat{h}_d & i = i' \\
\hat{h}_d & i \neq i',
\end{cases}$$
making the bottom level an LCP with a positive semidefinite coefficient matrix $\begin{bmatrix} B_d & I \\ -I & 0 \end{bmatrix}$. This makes the problem setting suitable for the method described in Section \ref{sec:posit-semid-m-x-dep}. but we will continue by presenting the strong duality approach to this problem.

\subsection{Strong duality reformulation of the trilevel electricity market model}
\label{app:mpc-pdc-model}

Using the primal-dual conversion rules for quadratic programs summarised in \textcite{dorn1960duality}, the dual of the bottom-level problem \eqref{eq:bottom-level-primal-arbitrage-reduced} can be stated as
\begin{subequations}
\label{eq:bottom-level-dual}
\begin{align}
  \min_{\lambda_{ijkd}, z_{ijkd}} \quad & \hat{h}_d Z_{id}^2 + \sum_{j \in J, k \in K} z_{ijkd}^{max} \lambda_{ijkd} \label{eq:bottom-level-dual.obj}\\
  \st \quad & - \lambda_{ijkd} \le \hat{h}d(Z_d+Z_{id}) - \sum_{k \in K} h_{kd} (\alpha_{kd}-w_{kd}) \nonumber \\
  &\qquad \qquad + (\gamma_{ijk} - w_{kd}) \ (z_{ijkd}) \ \forall j \in J, k \in K \label{eq:bottom-level-dual.con1}\\ 
  & z_{ijkd}, \lambda_{ijkd} \geq 0. \label{eq:bottom-level-dual.vars}
\end{align}
\end{subequations}

As described in Section \ref{sec:mpc-pdc}, we impose a \emph{strong duality} constraint stating that the objective value of the dual (minimisation) problem is less or equal to that of the primal (maximisation) problem, and combine constraints \eqref{eq:bottom-level-primal-arbitrage-reduced.con1}-\eqref{eq:bottom-level-primal-arbitrage-reduced.vars}, \eqref{eq:bottom-level-dual.con1}-\eqref{eq:bottom-level-dual.vars} and the strong duality constraint. A solution that satisfies these constraints must be optimal to \eqref{eq:bottom-level-primal-arbitrage-reduced} and \eqref{eq:bottom-level-dual}. Notice that the inequality version of the strong duality constraint is convex (as opposed to an equality constraint between the primal and dual objective values), and the other constraints are affine. 

Finally, we can write the primal and dual constraints and the strong duality constraint as
\begin{subequations}
\label{eq:bottom-level-pdc-arbitrage-reduced}
\begin{align}
  & z_{ijkd} \le z_{ijkd}^{max}  \quad (\lambda'_{ijkd})\ \forall i \in I, j \in J, k \in K \label{eq:bottom-level-pdc-arbitrage-reduced-con1} \\
  & - \lambda_{ijkd} - \hat{h}_d\left(Z_d + Z_{id}\right) \le - \sum_{k \in K}(\alpha_{kd}-w_{kd})h_{kd} \nonumber \\
  & \qquad \qquad + (\gamma_{ijkd} - w_{kd})  \quad (z'_{ijkd})\ \forall i \in I, j \in J, k \in K \label{eq:bottom-level-pdc-arbitrage-reduced-con2} \\
  & \sum_{i \in I} \left(\hat{h}_d Z_{id}^2 + \hat{h}_d Z_d Z_{id} + \sum_{j \in J, k \in K} z_{ijkd}^{max} \lambda_{ijkd} - \left( \sum_{k \in K} (\alpha_{kd} - w_{kd})h_{kd} \right) Z_{id} \right. \nonumber \\
  & \qquad \qquad \left.  + \sum_{j \in J, k \in K} \left(\gamma_{ijk}-w_{kd}\right) z_{ijkd}\right) \le 0 \quad (\epsilon_d)\ \label{eq:strong-duality-combined} \\
  & z_{ijkd}, \lambda_{ijkd} \geq 0, \label{eq:bottom-level-pdc-arbitrage-reduced.vars}
\end{align}
\end{subequations}
where the strong duality constraints for all producers $i \in I$ have been combined into a single constraint \eqref{eq:strong-duality-combined} to reduce the number of constraints as suggested in \textcite{pineda2018efficiently}.



The KKT conditions of the ISO problem \eqref{eq:profit_iso} combined with the constraints \eqref{eq:bottom-level-pdc-arbitrage-reduced} and the market-clearing constraint \eqref{eq:market_clearing_arb} are 
\begin{subequations}
\label{eq:iso_kkt_arb}
\begin{align}
  0 \le z_{ijkd} \perp &\lambda'_{ijkd} - \hat{h}_d(Z_d + Z_{id}) + \left(2\hat{h}_d\left(Z_d + Z_{id}\right) - \sum_{k \in K}(\alpha_{kd}-w_{kd})h_{kd} + \right. \nonumber \\
  & \left. \phantom{\sum_{k \in K}} (\gamma_{ijkd} - w_{kd}) \right) \epsilon_d -  (\mathbb{I}(j \in R) - \rho)\psi_{kd} \ge 0 \ \forall i \in I, j \in J, k \in K \\
  0 \le \lambda_{ijkd} \perp &-z'_{ijkd} + z_{ijkd}^{max}\epsilon_d \ge 0 \ \forall i \in I, j \in J, k \in K \\
  0 \le z'_{ijkd} \perp & \lambda_{ijkd} + \hat{h}_d(Z_d + Z_{id}) - \sum_{k \in K}(\alpha_{kd}-w_{kd})h_{kd} + (\gamma_{ijkd} - w_{kd}) \ge 0 \ \forall i \in I, j \in J, k \in K \\
  0 \le \lambda'_{ijkd} \perp & z_{ijkd}^{max} - z_{ijkd} \ge 0 \ \forall i \in I, j \in J, k \in K \\
  0 \le \phi_{ad}^- \perp & T_a^- + \sum_{k \in K} PTDF_{ka} y_{kd} \ge 0 \ \forall a \in A \\
  0 \le \phi_{ad}^+ \perp & T_a^+ - \sum_{k \in K} PTDF_{ka} y_{kd} \ge 0 \ \forall a \in A \\
  0 \le \psi_{kd} \perp & \sum_{i \in I, j \in R} z_{ijkd} - \rho\sum_{i \in I, j \in J} z_{ijkd} \ge 0 \ \forall k \in K \\
  0 \le \epsilon_d \perp & - \sum_{i \in I} \left(\hat{h}_d Z_{id}^2 + \hat{h}_d Z_D Z_{id} + \sum_{j \in J, k \in K} z_{ijkd}^{max} \lambda_{ijkd} - \right. \nonumber \\ 
  & \left. \left( \sum_{k \in K} (\alpha_{kd} - w_{kd})h_{kd}  \right) Z_{id} + \sum_{j \in J, k \in K} \left(\gamma_{ijk}-w_{kd}\right) z_{ijkd}\right) \ge 0 \label{eq:strongduality-comp} \\
  w_{kd} = &\sum_{a \in A}PTDF_{ka} (\phi_{ad}^+ - \phi_{ad}^-) \ \forall k \in K \\
  y_{kd} = &- Z_{kd} + h_{kd} Z_d + \sum_{k' \in K} L^d_{k,k'}(\alpha_{k'd} - w_{k'd}) \ \forall k \in K.
\end{align}
\end{subequations}
The indicator term $\mathbb{I}(j \in R)$ is 1 if $j \in R$, 0 otherwise, and the variables $z'_{ijkd}$ and $\lambda'_{ijkd}$ are the dual variables of the primal and dual constraints from the producer level for the system operator problem, and $\epsilon$ is the dual variable for the strong duality constraint. 

\end{document}